\documentclass[letterpaper, 10pt, conference]{ieeeconf}

\IEEEoverridecommandlockouts                              
\usepackage{amsmath}
\usepackage{amssymb}
\usepackage{epsfig}
\usepackage{graphicx}
\usepackage{color}
\usepackage{psfrag}
\newtheorem{thm}{Theorem}

\newtheorem{prp}{Proposition}

\newtheorem{rmq}{Remark}

\title{\LARGE \bf
Distributed source identification for wave equations:\\ an observer-based approach
}

\author{Marianne Chapouly and Mazyar Mirrahimi
\thanks{This work was supported by the ``Agence Nationale de la Recherche'' (ANR), Projet Jeunes Chercheurs EPOQ2 number 4136.}
\thanks{M. Chapouly and M. Mirrahimi are with Projet SISYPHE, INRIA Paris-Rocquencourt,
        Domaine de Voluceau, Rocquencourt B.P. 105, 78153 Le Chesnay Cedex, France
        {\tt\small marianne.chapouly@inria.fr} and {\tt\small mazyar.mirrahimi@inria.fr}}%
}

\begin{document}

\maketitle
\thispagestyle{empty}
\pagestyle{empty}

\begin{abstract}

In this paper, we consider a wave equation on a bounded interval where the initial
conditions are known (are zero) and we are rather interested in identifying an unknown
source term $q(x)$ thanks to the measurement output $y$ which is the
Neumann derivative on one of the boundaries. We use a back and forth iterative procedure and construct well-chosen observers
which allow to retrieve $q$ from $y$ in the minimal observation time.
\end{abstract}


\section{INTRODUCTION}\label{sec:intro}
\indent In a recent work~\cite{c1}, Blum and Auroux proposed a new inversion algorithm for identifying the initial state of an observable system, based on the application of back-and-forth observers. Noting that we have only access to a measurement output on a fixed time interval $(0,T)$, the idea consists in proposing a first asymptotic observer for the system that will be applied in this time interval and a second one that will be applied to the system where the direction of time is reversed. These observers are then used iteratively to get a better estimate of the initial state after each back-and-forth iteration. If the two observer gains  are well-chosen, so that the whole back-and-forth procedure induces a contracting error dynamic, one can ensure the convergence of the estimator towards the initial state.\\  \indent In particular, Ramdani-et-al~\cite{c3,c4} have considered (for the case of wave equations) the theoretical study of this problem applying techniques borrowed from semigroups theory. \\

Here, we consider a similar problem to~\cite{c4} for a wave equation, where the initial conditions are known (are zero) and we are rather interested in identifying an unknown source term $q(x)$. Let also $T>0$,  $\omega\in\mathbb{R}$ and let $q\in H^2(0,1)\cap H^1_0(0,1)$. We consider the following system
\begin{equation}
\label{eq1}
\left\{
\begin{array}{l}
u_{tt}-u_{xx}=q(x)\cos(\omega t),\,(t,x)\in(0,T)\times(0,1),\\
u(t,0)=u(t,1)=0,\,t\in(0,T),\\
u(0,x)=u_t(0,x)=0,\,x\in(0,1),\\
y(t)=u_x(t,0),\,t\in(0,T),
\end{array}
\right. 
\end{equation}
where $(u,u_t)$ represents the state of the system and $y$ is the output. The term $cos (\omega t)q(x)$, where $\omega$ is a fixed (known) frequency and $q(x)$ is unknown, is considered to be an external force which varies harmonically.\\
\indent For any $q\in H^2(0,1)\cap H^1_0(0,1)$ there exists a unique solution $u\in C^1([0,T];H^1_0(0,1))\cap C^2([0,T];L^2(0,1))$ to \eqref{eq1} and $u_x(t,0)\in H^1(0,T)$ (see \cite[Remark 2]{c5}). It is moreover well-known from the work \cite{c5} of Yamamoto that this problem is well-posed in the sense that one can retrieve the source term $q(x)$ from the measurement of $y$ on the time interval $(0,T)$ if $T$ is large enough.  
Our aim here is to propose well-chosen observers which allow, using back and forth procedure, to retrieve $q$ from $y$ in the minimal observation time. More precisely, we prove the following result

\textbf{Claim 1}
We can construct efficient observers for which the back-and-forth algorithm is convergent and which allow, using the measurement output $y(t)$ over the time-interval $(0,2)$, to reconstruct the unknown source term $q(x)$.\\
\begin{rmq}\label{rmq1}
Let us point out that since the spatial domain is given by $(0,1)$, the minimum observability time is given by $T=2$.
\end{rmq}
 Note that, whenever the whole initial state $(u(0,.), u_t(0,.),q)$ is unknown, system \eqref{eq1} is not observable. In order to realize this, one can consider the simpler case where the source term $q$ is given by only the two first modes of the wave equation, $q(x)=q_1\sin(\pi x)+q_2\sin(2\pi x)$, where $q_1$ and $q_2$ are the unknown scalars. In this case, \eqref{eq1} becomes equivalent to two independent oscillators with different frequencies and with two unknown source terms $q_1\cos(\omega t)$ and $q_2\cos(\omega t)$. Moreover, the output is given by a linear combination of the position of the oscillators. This six dimensional system with one output is not observable. However, if we know the initial state of the oscillators, the two parameters $q_1$ and $q_2$ become identifiable.\\
The back-and-forth estimator allows us to take into account this knowledge of the initial state $(u(0,.), u_t(0,.))=(0,0)$ of the wave equation. On the contrary, if we had only used a forward observer, we would have lost the information on the initial state of the system and therefore there would have been no reason for the observer to converge to the real parameters.\\
 We prove the convergence of the algorithm using Lyapunov techniques and LaSalle's principle. One of the main difficulties comes from the fact that the precompactness of the trajectories is not ensured since we deal with an infinite dimensional system and we also have to use sharp mathematical estimates.\\

The paper is organized as follows.\\
In Section 2, we prove the equivalence between \eqref{eq1} and another system which consists in a system composed of a wave equation without any source term and an oscillator in cascade.  Here the unknown term to retrieve is the initial condition.
Section 3 is devoted to the proof of Claim 1 and is divided in two parts. In the first one we study the well-posedness of the considered systems. In the second one we prove the convergence of the proposed algorithm. 


\section{An equivalent estimation problem}
We begin with  introducing the following system,

\begin{equation}
\label{eq2}
\left\{
\begin{array}{l}
w_{tt}-w_{xx}=0,\,(t,x)\in(0,T)\times(0,1),\\
w(t,0)=w(t,1)=0,\,t\in(0,T),\\
w(0,x)=q(x),\,w_t(0,x)=0,\,x\in(0,1),\\
\end{array}
\right. 
\end{equation}
\begin{equation}
\label{eq2osc}
\left\{
\begin{array}{l}
\dot{z_1}(t)=z_2(t),\,t\in(0,T),\\
\dot{z_2}(t)=-\omega^2z_1(t)+w_x(t,0),\,t\in(0,T),\\
z_1(0)=y(0),\, z_2(0)=\dot{y}(0),\\
Y(t)=z_1(t),\,t\in(0,T),
\end{array}
\right. 
\end{equation}
where $(w,w_t,z_1,z_2)$ represents the state of the system and $Y$ is the output.  System \eqref{eq2}-\eqref{eq2osc} is nothing but a homogeneous wave equation coupled with an oscillator. The unknown $q$ is now the initial datum of the system.
We have the following results, whose proofs can be found in \cite{c1bis}\\
\begin{prp}\label{prop0}
There exists a unique solution $(w,w_t,z_1,z_2)$ to \eqref{eq2}-\eqref{eq2osc} with the following regularity $$w\in C([0,T];H^2(0,1)\cap H^1_0(0,1)),$$ $$w_t\in C([0,T];H^1_0(0,1)),$$ 
$$(z_1,z_2)\in H^2(0,T)\times H^1(0,T).$$ Moreover, $w$ satisfies the following hidden property
\begin{gather}
\label{extra}
w_x(t,0)\in L^2(0,T).
\end{gather}
\end{prp}
and \\
\begin{prp}
\label{prop1}
For any $T>0$, for any $q\in H^2(0,1)$, $$y= Y \text{ in } H^2(0,T),$$ where $y$ denotes the output of system \eqref{eq1} and $Y$ denotes the output of system \eqref{eq2}.
\end{prp}

As a consequence of Proposition~\ref{prop1}, we will now focus on system \eqref{eq2}-\eqref{eq2osc} and provide well-chosen observers for this system in order to prove Claim 1.


\section{Observer design}
As we have already mentioned it in the introduction, we are going to prove Claim 1 using back-and-forth observers. We first rewrite system \eqref{eq2}-\eqref{eq2osc} in the following way
\begin{equation}
\label{eq2*}
\left\{
\begin{array}{l}
w^1_{t}=w^2,\,(t,x)\in(0,T)\times(0,1),\\
w^2_{t}=w^1_{xx},\,(t,x)\in(0,T)\times(0,1),\\
w^1(t,0)=w^1(t,1)=0,\,t\in(0,T),\\
w^1(0,x)=q(x),\,w^2(0,x)=0,\,x\in(0,1),\\
\end{array}
\right. 
\end{equation}
\begin{equation}
\label{eq2osc*}
\left\{
\begin{array}{l}
\dot{z_1}(t)=z_2(t),\,t\in(0,T),\\
\dot{z_2}(t)=-\omega^2z_1(t)+w^1_x(t,0),\,t\in(0,T),\\
z_1(0)=y(0),\, z_2(0)=\dot{y}(0),\\
Y(t)=z_1(t),\,t\in(0,T),
\end{array}
\right. 
\end{equation}
where, from above results, $$w^1\in C([0,T];H^2(0,1)\cap H^1_0(0,1)),$$ $$w^2\in C([0,T];H^1_0(0,1))$$
and
$$(z^1,z^2)\in H^2(0,1)\times H^1(0,1).$$
Let us recall that our aim is to propose an algorithm which allows to retrieve the unknown $q$ from the measurement output $Y$ on the time interval $(0,T)$. Our idea consists in designing a well-chosen asymptotic observer, which ensures the decrease of a same Lyapunov function in the back-and-forth procedure. To this aim, we first make the previous system periodical. More precisely, we define $(W^1,W^2,Z^1,Z^2)$ on $((0,+\infty)\times(0,1))^2\times (0,T)^2$ as the solution  to the following periodical system
\begin{equation}
\label{eq2**}
\left\{
\begin{array}{l}
W^1_{t}=W^2,\,(t,x)\in(2kT,(2k+1)T)\times(0,1),\vspace{0,1cm}\\
W^2_{t}=W^1_{xx},\,(t,x)\in(2kT,(2k+1)T)\times(0,1),\vspace{0,1cm}\\
W^1_{t}=-W^2,\,(t,x)\in((2k+1)T,(2k+2)T)\times(0,1),\vspace{0,1cm}\\
W^2_{t}=-W^1_{xx},\,(t,x)\in((2k+1)T,(2k+2)T)\times(0,1),\vspace{0,1cm}\\
W^1(t,0)=W^1(t,1)=0,\,t\in(0,+\infty),\vspace{0,1cm}\\
W^1(0,x)=q(x),\,W^2(0,x)=0,\,x\in(0,1),
\end{array}
\right. 
\end{equation}
\begin{equation}
\label{eq2osc**}
\left\{
\begin{array}{l}
\dot{Z_1}(t)=Z_2(t),\,t\in(2kT,(2k+1)T),\vspace{0,1cm}\\
\dot{Z_2}(t)=-\omega^2Z_1(t)+W^1_x(t,0),\,t\in(2kT,(2k+1)T),\vspace{0,1cm}\\
\dot{Z_1}(t)=-Z_2(t),\,t\in((2k+1)T,(2k+2)T),\vspace{0,1cm}\\
\dot{Z_2}(t)=\omega^2Z_1(t)-W^1_x(t,0),\\
\phantom{ttttttlmmlllmmmmllltt}t\in((2k+1)T,(2k+2)T),\vspace{0,1cm}\\
Z_1(0)=y(0),\, Z_2(0)=\dot{y}(0),\vspace{0,1cm}\\
Y(t)=Z_1(t).
\end{array}
\right. 
\end{equation}
for  $k\in \mathbb{N}$.
Using previous results, one easily sees that $$W^1\in L^{\infty}_{loc}(\mathbb{R}_+;H^2(0,1)\cap H^1_0(0,1)),$$ $$W^2\in L^{\infty}_{loc}(\mathbb{R}_+;H^1_0(0,1)),$$ 
and $$(Z^1,Z^2)\in H^2_{loc}(\mathbb{R}_+)\times H^1_{loc}(\mathbb{R}_+).$$ 
Indeed, the above periodic system represents the system \eqref{eq2*}-\eqref{eq2osc*} on the time intervals $(2kT,(2k+1)T)$, and the same system in the time-reversed manner on the intervals $((2k+1)T,(2k+2)T)$. In particular, one has
\begin{equation*}
\left\{
\begin{array}{l}
Y(t)=Y(t-2kT),\,t\in(2kT,(2k+1)T),\vspace{0,1cm}\\
Y(t)=Y((2k+2)T-t),\,t\in((2k+1)T,(2k+2)T).
\end{array}
\right. 
\end{equation*}
The asymptotic observer we propose is the following
\small
\begin{equation}
\label{eq13}
\left\{
\begin{array}{l}
\hat{W}^1_{t}=\hat{W}^2,\,(t,x)\in(2kT,(2k+1)T)\times(0,1),\vspace{0,1cm}\\
\hat{W}^2_{t}=\hat{W}^1_{xx},\,(t,x)\in(2kT,(2k+1)T)\times(0,1),\vspace{0,1cm}\\
\hat{W}^1_{t}=-\hat{W}^2,\,(t,x)\in((2k+1)T,(2k+2)T)\times(0,1),\vspace{0,1cm}\\
\hat{W}^2_{t}=-\hat{W}^1_{xx},\,(t,x)\in((2k+1)T,(2k+2)T)\times(0,1),\vspace{0,1cm}\\
\hat{W}^{1}(t,0)=\gamma_1(\hat{Z}_1(t)-Y(t))\vspace{0,1cm}\\
\phantom{ttttttttttttttt}+\gamma_1\gamma_2(\hat{Z}_3(t)-\int_0^{t}Y(s)\,ds)\phantom{ttt}t\in(0,\infty),\vspace{0,1cm}\\
\hat{W}^1(t,1)=0,\,t\in(0,+\infty),\vspace{0,1cm}\\
\hat{W}^1(0,x)=\hat{W}^2(0,x)=0,\,x\in(0,1),
\end{array}
\right. 
\end{equation}
\begin{equation}
\label{eq14}
\left\{
\begin{array}{l}
\dot{\hat{Z}}_1(t)=\hat{Z}_2(t)-\gamma_2(\hat{Z}_1(t)-Y(t)),\,t\in(2kT,(2k+1)T),\vspace{0,1cm}\\
\dot{\hat{Z}}_2(t)=-\omega^2\hat{Z}_1(t)+\hat{W}^1_x(t,0),\,\phantom{ggg}t\in(2kT,(2k+1)T),\vspace{0,1cm}\\
\dot{\hat{Z}}_1(t)=-\hat{Z}_2(t)-\gamma_2(\hat{Z}_1(t)-Y(t)),\vspace{0,1cm}\\
\phantom{gggggggggggggggggggggggg}t\in((2k+1)T,(2k+2)T),\vspace{0,1cm}\\
\dot{\hat{Z}}_2(t)=\omega^2\hat{Z}_1(t)-\hat{W}^1_x(t,0),\vspace{0,1cm}\\
\phantom{gggggggggggggggggggggggg}t\in((2k+1)T,(2k+2)T),\vspace{0,1cm}\\
\dot{\hat{Z}}_3(t)=\hat{Z}_1(t),\,\phantom{ggggggggggg}t\in(0,\infty),\vspace{0,1cm}\\
\hat{Z}_1(0)=\hat{Z}_2(0)=\hat{Z}_3(0)=0.
\end{array}
\right. 
\end{equation}\normalsize
Here $\gamma_1$ and $\gamma_2$, the observer gains, are strictly positive constants to be fixed as the design parameters.
Before studying the well-posedness of system \eqref{eq13}-\eqref{eq14}, let us introduce the error  dynamics. The error term being defined as the difference between the observer and the observed system, error  equations are the following
\begin{equation}
\label{eq17}
\left\{
\begin{array}{l}
\tilde{W}^1_{t}=\tilde{W}^2,\,(t,x)\in(2kT,(2k+1)T)\times(0,1),\vspace{0,1cm}\\
\tilde{W}^2_{t}=\tilde{W}^1_{xx},\,(t,x)\in(2kT,(2k+1)T)\times(0,1),\vspace{0,1cm}\\
\tilde{W}^1_{t}=-\tilde{W}^2,\,(t,x)\in((2k+1)T,(2k+2)T)\times(0,1),\vspace{0,1cm}\\
\tilde{W}^2_{t}=-\tilde{W}^1_{xx},\,(t,x)\in((2k+1)T,(2k+2)T)\times(0,1),\vspace{0,1cm}\\
\tilde{W}^{1}(t,0)=\gamma_1\tilde{Z}_1(t)+\gamma_1\gamma_2\tilde{Z}_3(t),\,t\in(0,\infty),\vspace{0,1cm}\\
\tilde{W}^1(t,1)=0,\,t\in(0,+\infty),\vspace{0,1cm}\\
\tilde{W}^1(0,x)=-q(x),\,\tilde{W}^2(0,x)=0,\,x\in(0,1),
\end{array}
\right. 
\end{equation}
\begin{equation}
\label{eq18}
\left\{
\begin{array}{l}
\dot{\tilde{Z}}_1(t)=\tilde{Z}_2(t)-\gamma_2\tilde{Z}_1(t),\,t\in(2kT,(2k+1)T),\vspace{0,1cm}\\
\dot{\tilde{Z}}_2(t)=-\omega^2\tilde{Z}_1(t)+\tilde{W}^1_x(t,0),\,t\in(2kT,(2k+1)T),\vspace{0,1cm}\\
\dot{\tilde{Z}}_1(t)=-\tilde{Z}_2(t)-\gamma_2\tilde{Z}_1(t),\vspace{0,1cm}\\
\phantom{ttttttttttttttttttttttttt}t\in((2k+1)T,(2k+2)T),\vspace{0,1cm}\\
\dot{\tilde{Z}}_2(t)=\omega^2\tilde{Z}_1(t)-\tilde{W}^1_x(t,0),\vspace{0,1cm}\\
\phantom{ttttttttttttttttttttttttt}t\in((2k+1)T,(2k+2)T),\vspace{0,1cm}\\
\dot{\tilde{Z}}_3(t)=\tilde{Z}_1(t),\,t\in(0,\infty),\vspace{0,1cm}\\
\tilde{Z}_1(0)=-y(0),\,\tilde{Z}_2(0)=-\dot{y}(0),\,\tilde{Z}_3(0)=0.
\end{array}
\right. 
\end{equation}


\subsection{Well-posedness}
From now and until the end, we denote $$H^1_r(0,1):=\{v\in H^1(0,1) \text{ s.t. } v(1)=0\}.$$ We have the following proposition (see \cite{c1bis} for the proof).\\
\begin{prp}
\label{prop2}
For any $T>0$, for any $(q^0,q^1)\in H^2(0,1)\cap H^1_0(0,1)\times H^1_0(0,1)$ and any $(\xi_1^0,\xi_2^0,\xi_3^0)\in\mathbb{R}^3$, there exists a unique solution $(v_1,v_2,\xi_1,\xi_2,\xi_3)$  to the following periodical Cauchy problem
\begin{equation}
\label{eqerror}
\left\{
\begin{array}{l}
v^1_{t}=v^2,\,(t,x)\in(2kT,(2k+1)T)\times(0,1),\vspace{0,1cm}\\
v^2_{t}=v^1_{xx},\,(t,x)\in(2kT,(2k+1)T)\times(0,1),\vspace{0,1cm}\\
v^1_{t}=-v^2,\,(t,x)\in((2k+1)T,(2k+2)T)\times(0,1),\vspace{0,1cm}\\
v^2_{t}=-v^1_{xx},\,(t,x)\in((2k+1)T,(2k+2)T)\times(0,1),\vspace{0,1cm}\\
v^1(t,1)=0,\,v^1(t,0)=\gamma_1\xi_1(t)+\gamma_1\gamma_2\xi_3(t),\vspace{0,1cm}\\
\phantom{nnnnnnnnnnnnnnnnn}t\in((2k+1)T,(2k+2)T),\vspace{0,1cm}\\
v^1(0,x)=q^0(x),\,v^2(0,x)=q^1(x),\,x\in(0,1),
\end{array}
\right. 
\end{equation}
\begin{equation}
\label{eqerrorbis}
\left\{
\begin{array}{l}
\dot{\xi_1}(t)=\xi_2(t)-\gamma_2\xi_1(t),\,t\in(2kT,(2k+1)T),\vspace{0,1cm}\\
\dot{\xi_2}(t)=-\omega^2\xi_1(t)+v^1_x(t,0),\,t\in(2kT,(2k+1)T),\vspace{0,1cm}\\
\dot{\xi_1}(t)=-\xi_2(t)-\gamma_2\xi_1(t),\vspace{0,1cm}\\
\phantom{ttllllldtddddddddlllttt}t\in((2k+1)T,(2k+2)T),\vspace{0,1cm}\\
\dot{\xi_2}(t)=\omega^2\xi_1(t)-v^1_x(t,0),\vspace{0,1cm}\\
\phantom{ttllllldtddddddddlllttt}t\in((2k+1)T,(2k+2)T),\vspace{0,1cm}\\
\dot{\xi_3}(t)=\xi_1(t),\,t\in(0,\infty),\vspace{0,1cm}\\
\xi_1(0)=\xi_1^0,\, \xi_2(0)=\xi_2^0,\, \xi_2(0)=\xi_3^0,
\end{array}
\right. 
\end{equation}
where $k\in \mathbb{N}$, with the following regularity

 $$v^1\in C([0,+\infty);H^2(0,1)\cap H^1_r(0,1)),$$ $$v_2\in C([0,+\infty);H^1_r(0,1)),$$
$$(\xi_1,\xi_2,\xi_3)\in H^1([0,+\infty))\times L^2([0,+\infty))\times H^2([0,+\infty)).$$

Moreover $(v^1,v^2,\xi_1,\xi_2,\xi_3)$ satisfies the following energy identities: for any $t\in (0,+\infty)$,
\begin{eqnarray*}
\begin{array}{l}
\label{eq26}
\vert v^1_x(t,.)\vert^2_{L^2(0,1)}+\vert v^2(t,.)\vert^2_{L^2(0,1)}+\gamma_1\vert\xi_2(t)\vert^2\vspace{0,2cm}\\
\phantom{ttttttt}+\gamma_1\omega^2\vert\xi_1(t)\vert^2+2\gamma_1\gamma_2\omega^2\vert\xi_1\vert^2_{L^2(0,t)}=\vspace{0,2cm}\\
\phantom{ttttttdddt}\vert q^0_x\vert^2_{L^2(0,1)}+\vert q^1\vert^2_{L^2(0,1)}+\gamma_1\vert\xi_2^0\vert^2+\gamma_1\omega^2\vert\xi_1^0\vert^2.
\end{array}
\end{eqnarray*}

and \small
\begin{eqnarray*}
\begin{array}{l}
\label{eq26bis}
\displaystyle \frac{1}{2}(\vert v^2_{t}(t,.)\vert^2_{L^2(0,1)}+\vert v^2_{x}(t,.)\vert^2_{L^2(0,1)}+\gamma_1\omega^4\vert\xi_1(t)\vert^2)\vspace{0,2cm}\\
\displaystyle+\frac{1}{4}\gamma_1\vert v^1_x(t,0)\vert^2
+\frac{1}{2}\gamma_1\gamma_2\omega^2\int_0^t\vert \xi_2(s)\vert^2\,ds+\gamma_1\omega^2\vert \xi_2(t)\vert^2\le \\
\phantom{ttttdttddtttttttttt}C(\vert q^0\vert^2_{H^2(0,1)}+\vert q^1\vert^2_{H^1(0,1)}+\vert \xi_1^0\vert^2+\vert \xi_2^0\vert^2),
\end{array}
\end{eqnarray*}\normalsize
where $C$ denotes a positive constant which only depends on $\omega, \gamma_1,\gamma_2$.\\
\end{prp}

One easily sees that the well-posedness of system \eqref{eq17}-\eqref{eq18} in 
\begin{eqnarray*}
\begin{array}{l}
C([0,+\infty);H^2(0,1)\cap H^1_r(0,1))\times C([0,+\infty);H^1_r(0,1))\\
\phantom{tttttttt}\times H^1([0,+\infty))\times L^2([0,+\infty))\times H^2([0,+\infty))
\end{array}
\end{eqnarray*}
directly follows from Proposition~\ref{prop2}. The well-posedness of system \eqref{eq13}-\eqref{eq14} in 
\begin{eqnarray*}
\begin{array}{l}
L^{\infty}_{loc}(\mathbb{R}_+;H^2(0,1)\cap H^1_r(0,1))\times L^{\infty}_{loc}(\mathbb{R}_+;H^1_r(0,1))\\
\phantom{tttttttt}\times H^1_{loc}(\mathbb{R}_+)\times L^2_{loc}(\mathbb{R}_+)\times H^2_{loc}(\mathbb{R}_+)
\end{array}
\end{eqnarray*}
 then follows (using also the well-posedness of system \eqref{eq2**}-\eqref{eq2osc**} in the same space).


\subsection{Asymptotic analysis}
The main goal of this section is to prove the following result, implying Claim 1\\
\begin{thm}\label{thm:conv}
For any $T\ge 2$,
\begin{gather}
\label{th2eq1}
\lim\limits_{n\rightarrow +\infty}\tilde{W}^1(2nT,.)=0\text{ in } H^1_r(0,1),\\
\label{th2eq1bis}
\lim\limits_{n\rightarrow +\infty} \tilde{W}^2(2nT,.) =0\text{ in }  L^2(0,1),\\
\label{th2eq2}
\lim\limits_{n\rightarrow +\infty}(Z_1(2nT),Z_2(2nT),Z_3(2nT))= (0,0,0).
\end{gather}
\end{thm}
\vspace{.2cm}

\textit{Proof of Theorem~\ref{thm:conv}.} Let us assume that \eqref{th2eq1}-\eqref{th2eq2} do not hold. Then, there exist a positive constant $\alpha$ and  a subsequence $(\phi(n))_{n\ge0}\in \mathbb{N}^\mathbb{N}$  with $\lim\limits_{n\rightarrow +\infty}\phi(n)=+\infty$ such that for any $n\ge0$, \small

\begin{multline}\label{th2eq3}
\vert (\tilde{W}^1(2\phi(n)T,.),\tilde{W}^2(2\phi(n)T,.))\vert_{H^1_r(0,1)\times L^2(0,1)} +\\
\vert (Z_1(2\phi(n)T),Z_2(2\phi(n)T),Z_3(2\phi(n)T))\vert >\alpha.
\end{multline}\normalsize

From Proposition~\ref{prop2}, $$\tilde{W}^1\in L^{\infty}(\mathbb{R}_+;H^2(0,1)\cap H^1_r(0,1)),$$ $$\tilde{W}^2\in L^{\infty}(\mathbb{R}_+;H^1_r(0,1)),$$ 
$$(\tilde{Z}_1,\tilde{Z}_2,\tilde{Z}_3)\in H^1(\mathbb{R}_+)\times L^2(\mathbb{R}_+)\times H^2(\mathbb{R}_+),$$

and for any time $t\in \mathbb{R}_+^*$, 
\begin{eqnarray}
\begin{array}{l}
\label{eq70}
\vert \tilde{W}^1_x(t,.)\vert^2_{L^2(0,1)}+\vert \tilde{W}^2(t,.)\vert^2_{L^2(0,1)}+\gamma_1\vert\tilde{Z}_2(t)\vert^2\vspace{0,2cm}\\
+\gamma_1\omega^2\vert\tilde{Z}_1(t)\vert^2+2\gamma_1\gamma_2\omega^2\vert\tilde{Z}_1\vert^2_{L^2(0,t)}=\vspace{0,2cm}\\
\phantom{ttttttttttttttttlt}\vert q\vert^2_{H^1(0,1)}+\gamma_1\vert \dot{y}(0)\vert^2+\gamma_1\omega^2\vert y(0)\vert^2.
\end{array}
\end{eqnarray}
Since $\{(\tilde{W}^1(2\phi(n)T,.),\tilde{W}^2(2\phi(n)T,.)),\,n\in\mathbb{N}\}$ is bounded in $H^2(0,1)\times H^1_r(0,1)$, it follows from Kato-Rellich's theorem that there exists a subsequence of $(\phi(n))_{n\ge0}$, that, for convenience, we still denote $(\phi(n))_{n\ge0}$, and there exists $(W^{\infty,1},W^{\infty,2})\in H^1_r(0,1)\times L^2(0,1)$ such that 
\begin{gather}
\label{cv}
\lim\limits_{n\rightarrow +\infty}\tilde{W}^1(2\phi(n)T,.)=W^{\infty,1} \text{ in } H^1_r(0,1),\\
\label{cv1}
\lim\limits_{n\rightarrow +\infty} \tilde{W}^2(2\phi(n)T,.)=W^{\infty,2} \text{ in }  L^2(0,1).
\end{gather}
From \eqref{eq70},  $\tilde{W}^1(2\phi(n)T,0)$, $\tilde{Z}_1(2\phi(n)T)$ and $\tilde{Z}_2(2\phi(n)T)$ are bounded, and from the fifth equation in \eqref{eq17} so is $\tilde{Z}_3(2\phi(n)T)$. Thus, (up to a subsequence of $(\phi(n))_{n\ge0}$),  there exists $(Z_1^{\infty},Z_2^{\infty},Z_3^{\infty})\in \mathbb{R}^3$ such that
\begin{gather}
\label{cvbis}
\lim\limits_{n\rightarrow +\infty} \tilde{Z}_1(2\phi(n)T)=Z_1^{\infty},\\
\label{cvbis1}
\lim\limits_{n\rightarrow +\infty} \tilde{Z}_2(2\phi(n)T)=Z_2^{\infty},\\
\label{cvbis2}
\lim\limits_{n\rightarrow +\infty} \tilde{Z}_3(2\phi(n)T))=Z_3^{\infty}.
\end{gather}

Let now $$(v^1,v^2)\in C([0,T];H^1_r(0,1))\times C([0,T];L^2(0,1)),$$ $$(x_1,x_2,x_3)\in L^2(0,1)\times L^2(0,1)\times L^2(0,1)$$ be solution to
\begin{equation}
\label{eq80}
\left\{
\begin{array}{l}
v^1_{t}=v^2,\,(t,x)\in(0,T)\times(0,1),\\
v^2_{t}=v^1_{xx},\,(t,x)\in(0,T)\times(0,1),\\
v^1(t,1)=0,\,v^1(t,0)=\gamma_1x_1(t)+\gamma_1\gamma_2x_3(t),\,t\in(0,T),\\
\dot{x_1}(t)=x_2(t)-\gamma_2x_1(t),\,t\in(0,T),\\
\dot{x}_2(t)=-\omega^2x_1(t)+v^1_x(t,0),\,t\in(0,T),\\
\dot{x}_3(t)=x_1(t),\,t\in(0,T),\\
v^1(0,x)=W^{\infty,1}(x),\,v^2(0,x)=W^{\infty,2}(x),\,x\in(0,1),\\
x_1(0)=Z_1^{\infty},\,x_2(0)=Z_2^{\infty},\,x_3(0)=Z_3^{\infty}.
\end{array}
\right. 
\end{equation}
Finally, let us define, for any $n\ge0$, $$(v^{1,n},v^{2,n})\in C([0,T];H^1_r(0,1))\times C([0,T];L^2(0,1))$$ and $$(x_1^n,x_2^n,x_3^n)\in L^2(0,1)\times L^2(0,1)\times L^2(0,1)$$ by
\begin{gather}
 \label{eq86}
v^{1,n}(t,.):=(\tilde{W}^1(2\phi(n)T+t,.),\\
 \label{eq86bis}
v^{2,n}(t,.):=\tilde{W}^2(2\phi(n)T+t,.)
\end{gather}
and 
\begin{gather}
 \label{eq87}
x_1^n(t):=\tilde{Z}_1(2\phi(n)T+t),\\
\label{eq87bis}
x_2^n(t):=\tilde{Z}_2(2\phi(n)T+t),\\
\label{eq87ter}
x_3^n(t)):=\tilde{Z}_3(2\phi(n)T+t)
\end{gather}
 for any $t\in(0,T)$. With such a definition, for any $n\ge0$, $(v^{1,n},v^{2,n},x_1^n,x_2^n,x_3^n)$ is solution to
\begin{equation}
\label{eq81}
\left\{
\begin{array}{l}
v^{1,n}_{t}=v^{2,n},\,(t,x)\in(0,T)\times(0,1),\vspace{0,1cm}\\
v^{2,n}_{t}=v^{1,n}_{xx},\,(t,x)\in(0,T)\times(0,1),\vspace{0,1cm}\\
v^{1,n}(t,1)=0,\,v^{1,n}(t,0)=\gamma_1x^n_1(t)+\gamma_1\gamma_2 x^n_3(t),\\
\phantom{hhhhhhhh}t\in(0,T),\vspace{0,1cm}\\
\dot{x}_1^n(t)=x_2^n(t)-\gamma_2 x_1^n(t),\,t\in(0,T),\vspace{0,1cm}\\
\dot{x}^n_2(t)=-\omega^2x_1^n(t)+v^{1,n}_x(t,0),\,t\in(0,T),\vspace{0,1cm}\\
\dot{x}^n_3(t)=x_1^n(t),\,t\in(0,T),\vspace{0,1cm}\\
v^{1,n}(0,x)=\tilde{W}^1(2\phi(n)T,x),\,x\in(0,1),\\
v^{2,n}(0,x)=\tilde{W}^2(2\phi(n)T,x),\,x\in(0,1),\vspace{0,1cm}\\
x_1^n(0)=\tilde{Z}_1(2\phi(n)T),\vspace{0,1cm}\\
x_2^n(0)=\tilde{Z}_2(2\phi(n)T),\vspace{0,1cm}\\
x_3^n(0)=\tilde{Z}_3(\phi(n)T).
\end{array}
\right. 
\end{equation}
From \eqref{eq70}-\eqref{eq80} and \eqref{eq81} (more precisely by continuity of flow with respect to the initial state),
\begin{gather}
 \label{eq84}
\lim\limits_{n\rightarrow +\infty}v^{1,n}= v^1\text{ in } L^{\infty}((0,T);H^1_r(0,1)),\\
\label{eq84b}
\lim\limits_{n\rightarrow +\infty} v^{2,n}=v^2\text{ in } L^{\infty}((0,T);L^2(0,1)),\\
\label{eq84bis}
\lim\limits_{n\rightarrow +\infty}(x_1^n,x_2^n,x_3^n)=(x_1,x_2,x_3)\text{ in } L^{\infty}(0,T)^3.
\end{gather}

We now introduce the following Lyapunov function\small
\begin{multline}
 \label{eq82}
\mathcal{V}(\tilde{W}^1,\tilde{W}^2,\tilde{Z}_1,\tilde{Z}_2,\tilde{Z}_3)(t):=\\
\displaystyle\frac{1}{2}\Big(\int_0^1(\vert \tilde{W}^1_x(t,x)\vert^2+\vert \tilde{W}^2(t,x)\vert^2)\,dx
\\
+\gamma_1\omega^2\vert \tilde{Z}_1(t)\vert^2+\gamma_1\vert \tilde{Z}_2(t)\vert^2\Big),\,t\in(0,T).
\end{multline}\normalsize
Indeed, $\mathcal{V}(t)\ge0,\,t\ge0$ and using \eqref{eq17} and \eqref{eq18}, one can compute, for any time $t\in(0,T)$,
\begin{gather}
 \label{eq82bis}
\frac{d}{dt}\mathcal{V}(\tilde{W}^1,\tilde{W}^2,\tilde{Z}_1,\tilde{Z}_2,\tilde{Z}_3)(t)=-\gamma_1\gamma_2\omega^2\vert \tilde{Z}_1(t)\vert^2\le 0.
\end{gather}
The function $\mathcal{V}$ is positive, decreasing. Consequently there exists $l\ge0$ such that
\begin{gather}
\label{eq88}
\lim\limits_{t\rightarrow +\infty}\mathcal{V}(\tilde{W}^1,\tilde{W}^2,\tilde{Z}_1,\tilde{Z}_2,\tilde{Z}_3)(t)=l.
\end{gather}

On the other side,  from \eqref{eq86}, \eqref{eq87}, \eqref{eq84} and \eqref{eq82}, one has, for any $t\in(0,T)$,
\begin{multline}
 \label{eq85}
  \lim\limits_{n\rightarrow+\infty}\mathcal{V}(\tilde{W}^1,\tilde{W}^2,\tilde{W}^2,\tilde{Z}_1,\tilde{Z}_2,\tilde{Z}_3)(2\phi(n)T+t)=\\
  \displaystyle\frac{1}{2}\Big(\int_0^1(\vert v^1_x(t,x)\vert^2+\vert v^2(t,x)\vert^2)\,dx\\
+\gamma_1\omega^2\vert x_1(t)\vert^2+\gamma_1\vert x_2(t)\vert^2\Big),
\end{multline}

Thus it follows from \eqref{eq88} and \eqref{eq85} that for any time $t\in(0,T)$,
$$\mathcal{V}(v^1,v^2,x_1,x_2,x_3)(t)=l.$$
In other words, $t\mapsto\mathcal{V}(v^1,v^2,x_1,x_2,x_3)(t)$  is constant on $(0,T)$ and thus
\begin{multline}
 \label{eq89}
\frac{d}{dt}\mathcal{V}(v^1,v^2,x_1,x_2,x_3)=\\-\gamma_1\gamma_2\omega^2\vert x_1(t)\vert^2=0,\,t\in(0,T).
\end{multline}
We finally obtain
\begin{gather}
 \label{eq90}
x_1\equiv0 \text{ on } (0,T).
\end{gather}
Consequently, from the fourth  and sixth lines in \eqref{eq80}, we also get
\begin{gather}
 \label{eq91}
x_2\equiv0,~\quad x_3\equiv  Z_3^\infty \quad \text{ on } (0,T),
\end{gather}
where $Z_3^\infty$ is a real constant.
Then, using the third and fifth lines in \eqref{eq80} and, \eqref{eq90} and \eqref{eq91} we see that
\begin{gather}
 \label{eq93}
v^1_x(t,0)\equiv0,\quad v^1(t,0)=\gamma_1\gamma_2 Z_3^\infty \quad\text{ on } (0,T).
\end{gather}
Consequently, \eqref{eq80} reduces to
\begin{equation}
\label{eq80bis}
\left\{
\begin{array}{l}
v^1_{t}=v^2,\,(t,x)\in(0,T)\times(0,1),\\
v^2_{t}=v^1_{xx},\,(t,x)\in(0,T)\times(0,1),\\
v^1(t,1)=v^1_x(t,0)=0,\,t\in(0,T),\\
v^1(t,0)=\gamma_1\gamma_2 Z_3^\infty ,\,t\in(0,T),\\
v^1(0,x)=W^{\infty,1}(x),\,v^2(0,x)=W^{\infty,2}(x),\,x\in(0,1),\\
x_1(t)=x_2(t)=0,~ x_3(t)=Z_3^\infty,\,t\in(0,T).
\end{array}
\right. 
\end{equation}

We recall that if $W^{\infty,1}(x)=\sum\limits_{k=1}^{+\infty}a_k\cos(\frac{(2k+1)\pi x}{2})$ and $W^{\infty,2}(x)=\sum\limits_{k=1}^{+\infty}b_k\cos(\frac{(2k+1)\pi x}{2})$, then $v^1$ can be expressed as\small
\begin{multline*}
 v^1(t,x)=\sum_{k=1}^{\infty}\Big[a_k\cos\left(\frac{(2k+1)\pi t}{2}\right)+\\
 \frac{2b_k}{(2k+1)\pi }\sin\left(\frac{(2k+1)\pi t}{2}\right)\Big]\cos\left(\frac{(2k+1)\pi x}{2}\right),\\
 (t,x)\in(0,T)\times(0,1).
\end{multline*}
At $x=0$, we have
\begin{multline*}
v^1(t,0)=Z_3^{\infty}=\sum_{k=1}^{\infty}\Big[ a_k\cos\left(\frac{(2k+1)\pi t}{2}\right)\\
+\frac{2b_k}{(2k+1)\pi }\sin\left(\frac{(2k+1)\pi t}{2}\right)\Big],\,t\in(0,T).
\end{multline*}
Multiplying this last equality by $\cos\left(\frac{(2k+1)\pi t}{2}\right)$, for any $k\in\mathbb{N}$, and integrating on $(0,2)$ we obtain 
\begin{gather*}
 \label{balard}
a_k=0,\,k\in\mathbb{N},
\end{gather*}
and thus, 
\begin{gather}
\label{syrtes}
W^{\infty,1}\equiv 0\text{ on } (0,1).
\end{gather}

In particular, 
\begin{gather}
\label{syrtes2}
W^{\infty,1}(0)=\gamma_1\gamma_2 Z^{\infty}_3=0
\end{gather}
 and $$v^1_t(t,0)=\sum_{k=1}^{\infty}b_k\cos(\frac{(2k+1)\pi t}{2}),\,t\in(0,T).$$

 Parseval's identity then implies
 \begin{multline} \label{eq100}
 \int_0^T \vert v^1_t(t,0)\vert^2\,dt\ge \displaystyle\int_0^2\vert v^1_t(t,0)\vert^2\,dt=\\ \sum\limits_{k=1}^{\infty} \vert b_k\vert^2
 =\vert W^{\infty,2}\vert^2_{L^2(0,1)}.
 \end{multline}
 Thus, from \eqref{eq80bis}, we finally get
$$W^{\infty,1}\equiv W^{\infty,2}\equiv0 \text{ on } (0,1)$$ 
and $$Z^{\infty}_1=Z^{\infty}_2=Z^{\infty}_3=0.$$
This is a contradiction with \eqref{th2eq3} and finishes the proof of the Theorem \ref{thm:conv}. 


\section{Numerical simulations}
\begin{figure}[h]\label{fig1}
\includegraphics[width=.5\textwidth]{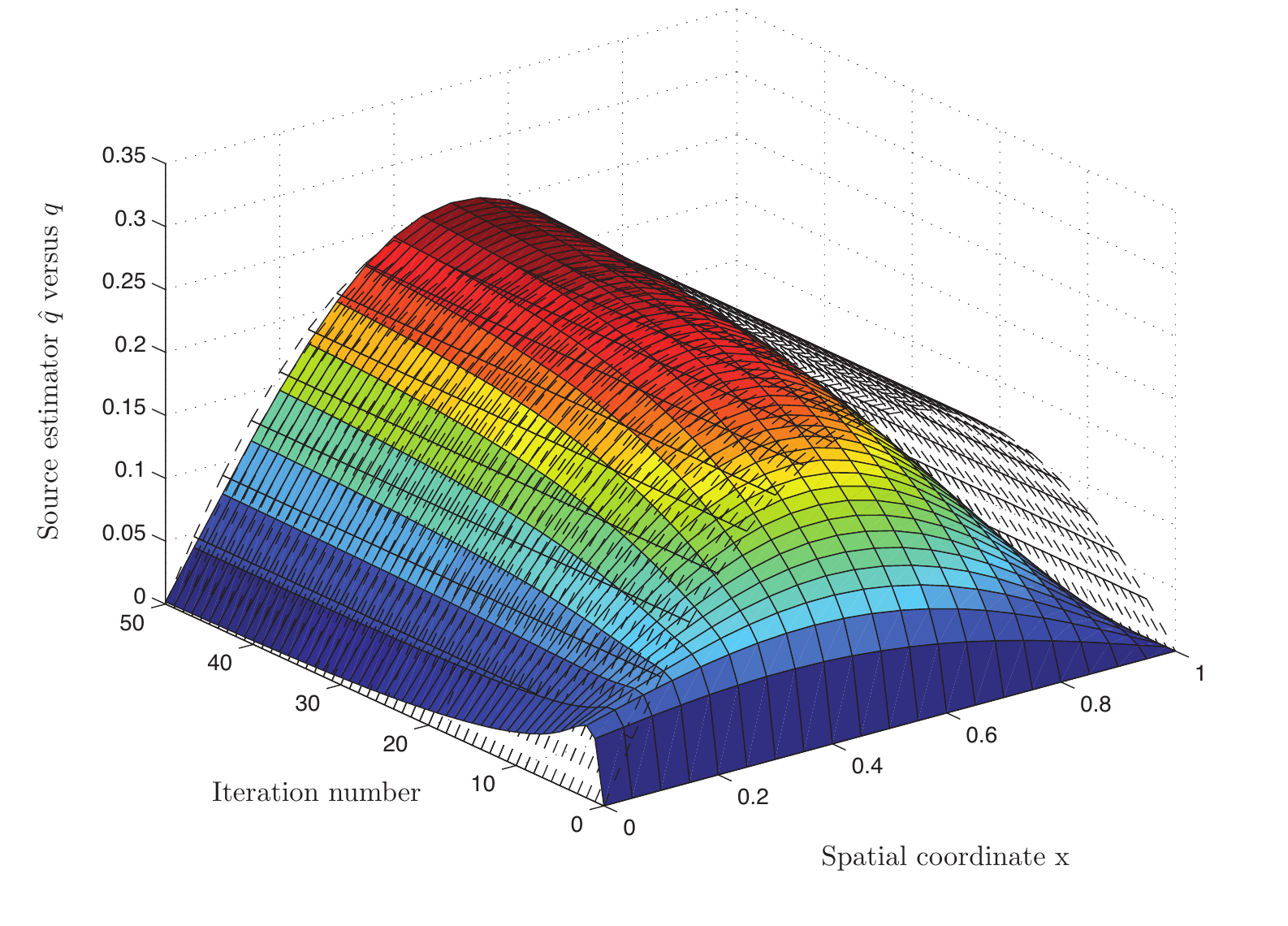}
\caption{The source estimator $\hat q$ is traced after each iteration; as it can be seen, the estimator after 50 iterations has converged towards $q$.}
\end{figure}
\begin{figure}[h]\label{fig2}
\includegraphics[width=.5\textwidth]{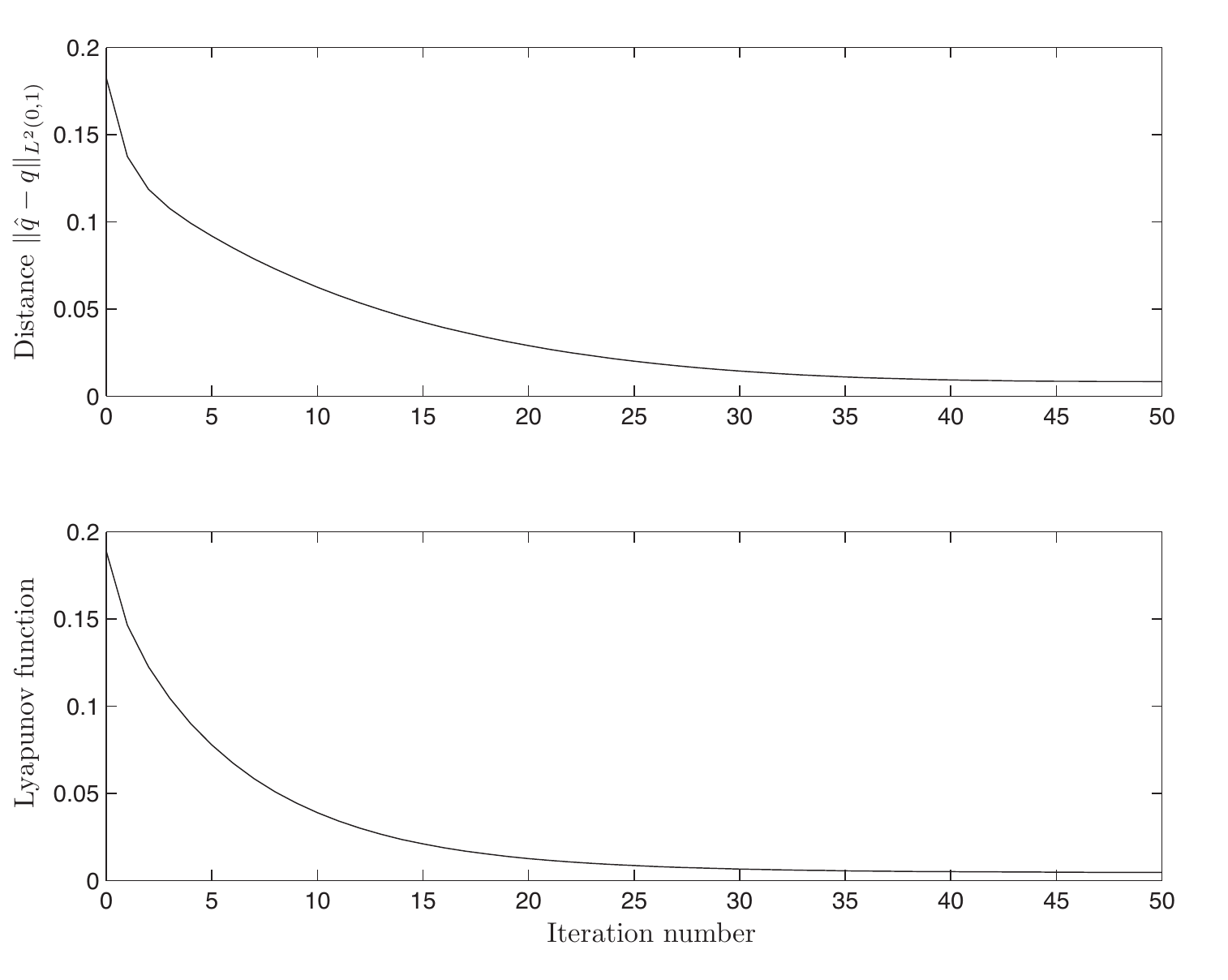}
\caption{The first plot illustrates the decrease of the $L^2$ distance between $\hat q$ and $q$ after each iteration and its convergence to zeros. The second plot illustrates the decrease of the Lyapunov function defined in~\eqref{eq82}.}
\end{figure}
In this section, we illustrate the efficiency of the above source estimation algorithm through numerical simulations. We consider the system~\eqref{eq1} with source term $q:=x-x^2$, together with the estimation algorithm~\eqref{eq13}-\eqref{eq14} with initial estimate $\hat{q}\equiv 0$. We fix the observation horizon to $T=3$ and consider 50 iterations of the estimator~\eqref{eq13}-\eqref{eq14}. The simulations of Figures~1and~2 illustrate the performance when we have added 10\% white noise on the measurement output and where the observer gains $\gamma_1$ and $\gamma_2$ are chosen to be
$$
\gamma_1=1,\qquad\gamma_2=1/2.
$$ 
The numerical simulations have been done through a finite difference method where the time and the space are discretized simultaneously. We have chosen a spatial discretization with 20 steps ($\Delta x=.05$) and a CFL coefficient of  $.005$ ($\Delta t=2.5e-04$).

\end{document}